\documentclass[10pt]{article} 
\usepackage[utf8]{inputenc} 
\usepackage{bbm}
\usepackage{relsize}
\usepackage{authblk}
\usepackage{geometry} 
 \usepackage{amsthm}
\usepackage{
	amsfonts,
	bbm,
	amsmath,
	appendix,
	color,
	amssymb,
	graphicx,
	bm,
	stmaryrd,
	mathabx,
	amssymb,
	mathdots,
	esint,
	cancel
}
\usepackage{amscd,pifont}
\usepackage{eucal,color,empheq}
\newcommand{\bea}{\begin{eqnarray}}
\newcommand{\eea}{\end{eqnarray}}
\newcommand{\R}{\mathbb{R}}
\newcommand{\E}{\mathbb{E}}

\newcommand{\pr}{\mathbbm{P}}
\newcommand{\lm}{\lambda_{max}}

\newcommand{\ind}{\mathbbm{1}}
\newcommand{\Z}{\mathbb{Z}}
\newtheorem{theorem}{Theorem}[section]

\newtheorem{lemma}[theorem]{Lemma}
\newtheorem{corollary}[theorem]{Corollary}

\newtheorem{example}{Example}[section]

\newtheorem{remark}[example]{Remark}

\title{Sharp asymptotics for Fredholm Pfaffians related to interacting particle systems and random matrices}
\author[1]{Will FitzGerald\thanks{will.fitzgerald@warwick.ac.uk}}
\author[1]{Roger Tribe\thanks{r.p.tribe@warwick.ac.uk}}
\author[1]{Oleg Zaboronski\thanks{olegz@maths.warwick.ac.uk}}
\affil[1]{\small{Department of Mathematics, 
University of Warwick, Coventry CV4~7AL, UK}}

\begin{document}
\maketitle

\abstract{It has been known since the pioneering paper of Mark Kac \cite{kac}, 
that the asymptotics of
Fredholm determinants can be studied using probabilistic
methods. We demonstrate the efficacy of Kac' approach by
studying the Fredholm Pfaffian describing the statistics 
of both non-Hermitian random matrices and annihilating Brownian motions.
Namely, we establish the following two results. 
Firstly, let $\sqrt{N}+\lambda_{max}$ be the largest real eigenvalue of a random $N\times N$ matrix
with independent $N(0,1)$ entries (the `real Ginibre matrix'). Consider
the limiting $N\rightarrow \infty$ distribution $\pr[\lm<-L]$ of the shifted maximal real eigenvalue $\lambda_{max}$. Then
\[
\lim_{L\rightarrow \infty}
e^{\frac{1}{2\sqrt{2\pi}}\zeta\left(\frac{3}{2}\right)L}
\pr\left(\lm<-L\right)
=e^{C_e},
\] 
where $\zeta$ is the Riemann zeta-function and 
\[
C_e=\frac{1}{2}\log 2+\frac{1}{4\pi}\sum_{n=1}^{\infty}\frac{1}{n}
\left(-\pi+\sum_{m=1}^{n-1}\frac{1}{\sqrt{m(n-m)}}\right).
\]
Secondly, let $X_t^{(max)}$ be the position of the rightmost particle at time $t$
for a system of annihilating Brownian motions (ABM's) started from every
point of $\R_{-}$. Then
\[
\lim_{L\rightarrow \infty}
e^{\frac{1}{2\sqrt{2\pi}}\zeta\left(\frac{3}{2}\right)L}
\pr\left(\frac{X_{t}^{(max)}}{\sqrt{4t}}<-L\right)
=e^{C_e}.
\] 
These statements are a sharp counterpart of the results of \cite{pop} improved by computing
the terms of order $L^{0}$ in the asymptotic expansion of the 
corresponding Fredholm Pfaffian.}

\section{Introduction and the main result}
The present paper continues the investigation of the statistics of the real eigenvalues
for random matrices with independent normal matrix elements (the so-called real
Ginibre ensemble) and particles for the system of annihilating Brownian
motions started in \cite{pop}. 

A mathematical way of describing
random arrangements of points representing eigenvalues or particle positions is the theory of point processes, see  \cite{daley} for a review.
Important subclasses of point processes are determinantal and Pfaffian point processes,
whose correlation functions are given by determinants or Pfaffians of kernels of certain
integral operators. Well known examples of determinantal point process are the laws of eigenvalues for random Hermitian, unitary and complex Gaussian matrix models;
the eigenvalues statistics for symmetric, symplectic and real random Gaussian matrices are described by Pfaffian point processes,
see \cite{zeitouni}, \cite{mehta} for reviews. Moreover, determinantal and Pfaffian point
processes describe the distribution of particles for a number of interacting
particle systems such as the totally asymmetric exclusion process  \cite{kurtj} and reaction-diffusion
systems for certain combinations of annihilation, coalescence, branching and immigration of particles,
\cite{bg_mp_rt_oz}, \cite{bg_rt_oz}, \cite{roger_oleg}. 

Of a particular importance for the current investigation is the fact that the law of the 
real eigenvalues for the real Ginibre ensemble is a Pfaffian point process, \cite{borodin_sinclair}, \cite{forrester_nagao},
\cite{sommers}. Moreover, its bulk scaling limit coincides (up to a diffusive rescaling)
with the fixed time law of annihilating Brownian motions started at every point of the real
line \cite{roger_oleg}, its edge scaling limit coincides with the fixed time law  of
annihilating Brownian motions started at every point of the negative part of the real line
\cite{bg_mp_rt_oz}, \cite{borodin_err}.

The probabilities of `gaps' (regions of space void of any particles) are a fundamental object for point
processes, which in fact characterise the law of a simple one-dimensional process uniquely. 
For determinantal and Pfaffian
point processes gap probabilities are given by Fredholm determinants and Pfaffians of integral 
operators determined by the kernels of the corresponding processes. A particular instance of   
gap probability is the distribution of the rightmost of leftmost particle for the process, 
for example the statsitics of the largest eigenvalue of a random 
Hermitian matrix (the Tracy-Widom distribution, \cite{mehta}), or the largest real eigenvalue of a real random matrix (the Rider-Sinclair distribution, \cite{rider_sinclair_2014}). 

 An exact calculation of a Fredholm determinant or a Pfaffian is impossible in all but a few 
 special cases. Fortunately, the asymptotics of gap probabilities in the limit of large empty intervals can be studied in many important cases. If, for example, the operator is translationally invariant, the asymptotics of the corresponding Fredholm determinant can be studied using Szeg\"o's theorem and its modifications, see \cite{szego} for review. In particular, Szeg\"o's
 theorem was used by Derrida and Zeitak to calculate the asymptotics of a single gap probabilities in coalescence-annihilation
 model started at every point of the real line \cite{derrida}. Note that for the purely
 annihilating case, Derrida-Zeitak's calculation is non-rigorous due to the presence
 of Hartwig-Fisher singularities in the kernel, but the final answer is believed to be correct
 and can be rigorised using an appropriate modification of Szeg\"o's formula, see
 \cite{szego}, Chapter 6. In \cite{forrester} Forrester used the Derrida-Zeitak formula
 and the connection between the real Ginibre random matrix model and 
annihilating Brownian motions stated above to calculate the asymptotics
of gap probabilities for the distribution of real eigenvalues in the bulk. For us it was a crucial
result which inspired our current research. 

 In the absence of 
 translational invariance, the situation is more complicated.
If the operator is integrable (as is the case for all point processes
describing eigenvalues of Hermitian random matrices),
the asymptotics of the distribution of extreme eigenvalues
 can be studied
by reducing the problem to a matrix Riemann-Hilbert problem
and analyzing the latter, see e.g. \cite{its}, Chapter XV for review.
As was discovered recently in \cite{bothner}, the operator $K$ which defines the Pfaffiant point
process for the annihilating
Brownian motions and the edge scaling limit of the real Ginibre ensemble
is conjugated to an integrable operator. This is a significant development, 
placing the real Ginibre ensemble firmly in the realm of integrable systems. 
Unfortunately, the associated Riemann-Hilbert problem turned out to be
rather complicated allowing the calculation of the asymptotic of the Fredholm
Pfaffian only for the operator $\gamma K$, where $\gamma <1$, thus making
it diffucult to relate the answer to the distribution of the largest real
eigenvalue ($\gamma=1$). Notice however, that in the context of particle system, $\gamma K$ for $\gamma<1$
has a clear probabilistic meaning - it 
describes the statistics of mixed annihilating-coalescing 
Brownian motions, \cite{bg_mp_rt_oz}.

An alternative approach to asymptotics of Fredholm determinants was pioneered
by Mark Kac \cite{kac} who was the first to state and prove a continuous version
of Szeg\"o's theorem, which originally was formulated for Toeplitz matrices rather than
translationally invariant operators. The main idea due to Kac is to interprete
the log-det expansion of the Fredholm determinant of the trace-class integral operator $T$ 
acting on $L^2$ functions on $I \subset \R$,
\[
\log \det(I-T)=-\sum_{n=1}^\infty \frac{1}{n} \int_{I^n} dx_1\ldots dx_n T(x_1- x_2)
T(x_2-x_3)\ldots T(x_{n}-x_1),
\]
 as a certain expectation with respect to the measure of a random walk whose
 increments have a (pseudo) distribution $T(x)dx$. Of course, the result of \cite{kac}
 can be derived directly by taking the continuous limit of Szeg\"o's 
 theorem for Toepltiz matrices.
 
 However, as we showed in \cite{pop},  
 the probabilistic approach can be used for the derivation of new results, namely
 the tails of the distribution of the rightmost real eigenvalue (the rightmost particle) 
 for the edge scaling limit of the real Ginibre ensemble (annihilating Brownian
 motions). This is the most difficult $\gamma=1$ case in the terminology of 
 \cite{bothner}, see the discussion above. In \cite{pop} we already calculated the asymptotic of the relevant Fredholm
 Pfaffian up to $O(1)$ errors.
 In the current paper we will show that by sticking closer to the 
 original Kac argument we can calculate the constant term as well as characterise
 the size of the correction. The calculation turns out to be rather short and intuitive. It is based some classical properties of random walks with general increments, as discussed in \cite{feller}.

The main result of the paper is the following statement. Let 
 $\sqrt{N}+\lambda_{max}$ be the largest real eigenvalue of the $N\times N$ real 
 Ginibre ensemble. Let $\E_N$ denote the ensemble expectation. Let
 \bea
 \pr(\lm<-L)=\lim_{N\rightarrow \infty} \E_{N}(\ind(\lm<-L))
 \eea
 be the edge scaling limit of the distribution of the largest real eigenvalue. 
\begin{theorem}\label{thm}
For $L>0$,
\bea\label{eq:tneg}
\lim_{L\rightarrow \infty}
e^{\frac{\zeta(3/2)}{2\sqrt{2\pi}}L}
\pr\left(\lm<-L\right)=e^{C_e},
\eea
where
\bea\label{edgeconst}
C_e=\frac{1}{2}\log 2+\frac{1}{4\pi}\sum_{n=1}^{\infty}\frac{1}{n}
\left(-\pi+\sum_{m=1}^{n-1}\frac{1}{\sqrt{m(n-m)}}\right).
\eea
More precisely,
\bea\label{logedge}
\log \pr\left(\lm<-L\right)=-\frac{\zeta(3/2)}{2\sqrt{2\pi}}L+C_e+o(L^{-1+}),
\eea
where for $o(L^{-1+})$: for any $\mu>0$, $\lim_{L\rightarrow \infty} L^{1-\mu}o(L^{-1+})=0$.
\end{theorem}
Let us analyse the presented asymptotic for the gap probability at the edge of the spectrum in more detail.
Numerically, $\exp(C_e)\approx 0.75$, which is consistent with its numerical value obtained in \cite{bothner}.
Let us also compare (\ref{logedge}) with the bulk scaling limit of the probabilty $\pr(N(-L,0)=0)$ that the interval $(-L,0)$
contains no real eigenvalues. As predicted by the Derrida-Zeitak formula \cite{derrida} applied to the real 
Ginibre ensemble in \cite{forrester},
\bea\label{logbulk}
\log \pr\left(N(-L,0)=0\right)=-\frac{\zeta(3/2)}{2\sqrt{2\pi}}L+C_b+o(L^{0}),
\eea
where
\bea\label{bulkconst}
C_b=\log 2+\frac{1}{4\pi}\sum_{n=1}^{\infty}\frac{1}{n}
\left(-\pi+\sum_{m=1}^{n-1}\frac{1}{\sqrt{m(n-m)}}\right).
\eea
Comparing (\ref{logbulk}) and (\ref{logedge}) we see that the leading terms coincide. This is not
very surprising, see \cite{pop} for a heuristic explanation. However, there is no reason why
the $O(1)$ terms should be the same. In fact, we see from the above formulae that
\bea\label{ratio}
\lim_{L\rightarrow \infty} \frac{\pr\left(\lm<-L\right)}{\pr\left(N(-L,0)=0\right)}=e^{C_e-C_b}=\frac{1}{\sqrt{2}}.
\eea
It would be interesting to see if it were possible to derive relation (\ref{ratio}) without computing
$\pr(\lm <-L)$ and $\pr(N(-L,0)=0)$ separately.

Interpreted in terms of particle systems, our result looks as follows:
\begin{corollary}\label{cor}
Consider the system of instantaneously annihilating Brownian motions on the real line started from every point of $\R_{-}$ (half-space maximal entrance law). 
Let $X_t^{(max)}$ be the position of the rightmost particle at a fixed time $t>0$. Then
\bea
\lim_{X\rightarrow \infty} e^{\frac{1}{2\sqrt{2\pi}}\zeta\left(\frac{3}{2}\right)\frac{X}{\sqrt{4t}}}\pr\left(X_t^{(max)}<-X\right)=e^{C_e}.
\eea
\end{corollary}
The Corollary is a direct consequence of the observation that 
\bea
X_t^{(max)}\stackrel{(d)}{\sim} \sqrt{4t} \lambda_{max}.
\eea
This in turn follows from the fact that the edge scaling limit of the law of real eigenvalues
for the real Ginibre ensemble and the single time distribution of  ABM's with half-space
maximal entrance law rescaled by $1/\sqrt{4t}$ can be characterised by the same Pfaffian
point process, see \cite{borodin_sinclair}, \cite{borodin_err} and \cite{bg_mp_rt_oz}
for the proof.

The rest of the paper is organised as follows. In Section \ref{proof-thm} we collect
the probabilistic tools necessary to establish our main result, explain the main idea
for the argument and finally prove Theorem \ref{thm}. In Section \ref{proof_lemmas}
we prove the probabilistic lemmas used to derive the statement of the Theorem. For the sake of
completeness we also present a streamlined proof of the key identity due to
Kac  \cite{kac}, which underpins our argument.
\section{The proof of Theorem \ref{thm}}\label{proof-thm}
The starting point for the proof is the Rider-Sinclair
formula \cite{rider_sinclair_2014}, which gives a Fredholm Pfaffian expression for $\pr(\lm<-L)$.
More specifically, we will use a probabilistic restatement of Rider-Sinclair's result
proved in \cite{pop}, which can be
explained as follows.
Let $\left(B_n,n\geq 0\right)$ be the discrete time random  walk with Gaussian 
$N(0,1/2)$ increments started at zero.
Let
\begin{equation} \label{eq:tau0}
\left\{ \begin{array}{rcl}
\tau_L &=&\inf_{n> 0} \{2n-1:~B_{2n-1}\geq L\},\\
\tau_0&=&\inf_{n> 0} \{2n:~B_{2n} \leq 0\}.
\end{array} \right. 
\end{equation}
In words: $\tau_0$ is the smallest $even$ time such that 
$B_{\tau_{0}}\leq 0$, $\tau_L$ is the smallest $odd$ time such that
$B_{\tau_L}\geq L$. 
Also, let 
\bea\label{eqsup}
M_{2n}=\sup\{B_k:k~is~odd,~k< 2n \}
\eea
be the infimum of the random walk $(B_k)_{k\geq 0}$ taken over all odd times not exceeding the time $2n$. Then
\begin{theorem} \label{thm:prob}
\bea\label{eq:probm}
\\
\nonumber
\pr(\lm<-L)=\sqrt{\pr(\tau_L<\tau_0)} \;
e^{-\frac{L}{2} \; \mathbb{E}\left(\delta_0\left(B_{\tau_0}\right)\right)}
e^{\frac{1}{2} \mathbb{E}\left(\min(L,M_{\tau_0})\delta_0\left(B_{\tau_0}\right)\right)}.
\eea
\end{theorem}
\begin{remark} We use the expression $\mathbb{E}(X\delta_{y}(Y))$ to mean 
a continuous Lebesgue density for the measure
$\mathbb{E}(X\mathbbm{1}(Y\in dy))$ evaluated at $y$. If $y=0$, we sometimes write 
$\mathbb{E}(X\delta_{0}(Y))$ as 
$\mathbb{E}(X\mathbbm{1}(Y\in d0))$.
\end{remark} 
\begin{remark}
Notice a slight change of notations in (\ref{eqsup}), (\ref{eq:probm}) in comparison with
formulae (1.7), (1.9) of \cite{pop}. 
\end{remark}
We will also need the following two facts from \cite{pop}:
\begin{lemma}\label{thm:exit}
As $L\rightarrow \infty$,
\bea\label{eq:mart}
\pr[\tau_L<\tau_0]= \frac{1}{\sqrt{2}L} (1+o(L^{-1+})).
\eea
Also,
\bea\label{eqlt}
\E\left(\delta_0\left(B_{\tau_0}\right)\right)=\frac{\zeta(3/2)}{\sqrt{2\pi}}.
\eea 
\end{lemma}
Formula (\ref{eq:mart}) is a slight improvement 
on Lemma 3.2 of \cite{pop}, which only claims the error bound
of magnitude $O(L^{-1/2})$. The improved bound is obtained
simply by using H\"older rather than the Cauchy-Schwarz inequality in the proof, without
changing the rest of the argument. Equation (\ref{eqlt})
is equally straightforward to check and we will do it below to illustrate
the utility of a probabilistic approach, see Remark \ref{leadtrm} below.
In addition, we need the following key statement due to Mark Kac \cite{kac}:
\begin{lemma}[Mark Kac, 1954]\label{lemma_kac}
Let $(X_i)_{i\geq 1}$ be independent identically distributed random variables having
continuous even density function $\rho$ on $\R$, and $S_k=X_1+X_2+\ldots+X_k, k\geq 1$.
Then
\bea 
\rho^{(n)}(0)\E\left( \max (0, S_1, S_2, \ldots, S_{n-1}) | S_n=0\right)=\frac{n}{2}\int_0^\infty
x \sum_{k=1}^{n-1} \frac{\rho^{(k)}(x)\rho^{(n-k)}(x)}{k(n-k)}dx,
\eea
where $\rho^{(k)}$ denotes the $k$-fold convolution of $\rho$ with itself or, in other words, the density
function of $S_k$.
\end{lemma}
Substituting (\ref{eq:mart}, \ref{eqlt}) into (\ref{eq:probm}), we find
\bea\label{meq1}
\log \pr(\lm<-L)+\frac{\zeta(3/2)}{2\sqrt{2\pi}}L=R(L),
\eea
where
\bea\label{rlterm}
R(L)=\frac{1}{2} \mathbb{E}\left(\min(L,M_{\tau_0})\delta_0\left(B_{\tau_0}\right)\right)
-\frac{1}{2}\log L-\frac{1}{4}\log 2+o(L^{-1+}),
\eea
where the notation $o(L^{-1/2+})$ means that for any $\mu>0$,
\bea
\lim_{L\rightarrow \infty} L^{1/2-\mu} o(L^{-1+})=0.
\eea
It remains to calculate the leading asymptotic of $R(L)$, which turns out to be $O(L^0)$.

Decomposing over the values of $\tau_0=2n$,
\bea
\E\left(\min(L,I_{\tau_0})\delta_0\left(B_{\tau_0}\right)\right)=\sum_{n=1}^\infty p_n(L),
\eea
where
\bea
p_n(L)=\E\left(\min(L,M_{2n})\ind(\tau_0=2n)\delta_{0}\left(B_{2n}\right) \right).
\eea
The summands $p_n(L), ~n\geq 1$, can be simplified using the cyclic 
invariance of the increments of the random walk $B$. Namely, we have the following result proved in Section  
\ref{proof_lemmas}.
\begin{lemma}\label{thm:cyclic}
Let $(S_n)_{n\geq 1}$ be a discrete time random walk such that the distribution of increments
has a continuous density.. Then
\bea\label{pncyc}
p_n(L)=\frac{1}{n}\E\left(\min\left(L, M_{2n}-m_{2n}\right)\delta_{0}(S_{2n})  \right),
\eea 
where 
\bea\label{eqinf}
m_{2n}=\inf\{S_k:k~is~even,~k\leq 2n \}
\eea
is the infimum of the random walk taken over even times not exceeding $2n$
and $M_n$ is the supremum of the walk over odd times defined in (\ref{eqsup}).
\end{lemma}
\begin{remark}
Notice that the above statement does not rely on the Gaussianity of increments.
It is a particular instance of a family of results for random walks conditioned to finish
at zero found in \cite{feller}.
\end{remark}
We conclude that $p_n$'s are fully determined by a joint distribution of the maximum, the minimum
and the final position of the random walk. 
Let us fix $\epsilon \in (0,2)$. Then
\bea\label{expect}
\E\left(\min(L,I_{\tau_0})\delta_0\left(B_{\tau_0}\right)\right)=
\sum_{n=1}^{\lfloor  L^{2-\epsilon} \rfloor} p_n(L)+
\sum_{\lfloor  L^{2-\epsilon} \rfloor+1}^\infty p_n(L).
\eea
For $n\leq L^{2-\epsilon}$, $M_{2n}-m_{2n}<L$ with probability close to $1$.
Therefore (\ref{pncyc}) is well approximated by $p_n(L) \approx \frac{1}{n}\E\left(\left(M_{2n}-m_{2n}\right)\delta_{0}(B_{2n})  \right)$, which can be computed adapting the original Kac' argument \cite{kac}. This approximation
decays as $1/2n$ at large $n$'s, leading to the logarithmic divergence of $\sum_{n=0}^\infty p_{n}(L)$
and thus the necessity for a separate analysis for large $n$.
Fortunately, for 
$n> L^{2-\epsilon}$, the random walk can be well approximated by a Brownian motion. Then
$p_n(L)$ can be computed using the classical Levy's result for the trivariate distribution of the 
supremum, the infimum and the final value of the Brownian motion on an interval, see e.g. 
\cite{billingsey}, \cite{feller} for review.  Rigorising the argument, we arrive at the following Lemma proved in Section \ref{proof_lemmas}.
\begin{lemma}\label{errbounds}
For $n\leq L^{2-\epsilon}$,
\bea\label{plow}
p_n(L)=\frac{1}{2\pi n} \sum_{k=1}^{n-1} \frac{1}{\sqrt{k(n-k)}}+E^{(1)}_n(L),
\eea
where 
\bea\label{e1}
|E^{(1)}_n(L)|\leq \sqrt{\frac{8}{\pi n^3}}
L \frac{e^{-\frac{2L^2}{n}}}{1-e^{-2L^{\epsilon}}}. 
\eea
For $n\geq L^{2-\epsilon}$, 
\bea\label{phi}
p_n(L)=\frac{1}{2n}-\sqrt{\frac{2}{\pi n^3}} L\sum_{k=1}^\infty e^{-\frac{2k^2L^2}{n}}+E^{(2)}_n(L),
\eea
where for any fixed $\gamma>0$, there exists an $n$-independent constant $C_{\gamma}>0$ such that 
\bea\label{e2}
|E^{(2)}_n(L)|\leq C_\gamma n^{-3/2+\gamma}.
\eea
\end{lemma}   
Substituting (\ref{plow}, \ref{phi}) into (\ref{expect}) and then into (\ref{rlterm}), we find
\bea\label{rlterm2}
R(L)&=&\frac{1}{4\pi}\sum_{n=1}^{\lfloor L^{2-\epsilon}\rfloor} \frac{1}{n}\sum_{k=1}^{n-1} \frac{1}{\sqrt{k(n-k)}}
+\frac{1}{2}\sum_{n=\lfloor L^{2-\epsilon} \rfloor+1}^\infty\left(\frac{1}{2n}-\sqrt{\frac{2}{\pi n^3}}L \Omega 
\left(\frac{2L^2}{\pi n} \right) \right)\nonumber\\
&-&\frac{1}{4}\log\left(2L^2\right)+
\frac{1}{2}\sum_{n=1}^{\lfloor L^{2-\epsilon} \rfloor} E_{n}^{(1)}(L)
+\frac{1}{2}\sum_{n=\lfloor L^{2-\epsilon} \rfloor+1}^\infty E_{n}^{(2)}(L)+o\left(L^{-1/2+}\right),
\eea 
where
\bea
\Omega(t):=\sum_{k=1}^\infty e^{-\pi k^2 t }
\eea
is a function on $\R$ closely related to Jacobi's $\theta$-function,
$\Omega(t)=\frac{\theta(0,it)-1}{2}$, see \cite{wwt} for a review.
 Recall that $\Omega(t) =O(e^{-\pi t})$ for $t\rightarrow \infty$
 and $\Omega(t)=O(t^{-1/2})$ for $t\downarrow 0$, which makes it
 easy to verify the convergence of integral bounds derived below.

First, let us estimate the error term in (\ref{rlterm2}) using (\ref{e1}, \ref{e2}).
Notice that the function $f(x)=x^{-3/2+\gamma}$ is decreasing on $\R_+$ and the function
$g(x)=x^{-3/2}\exp(-2k^2L^2/x)$ is increasing  for $0<x<L^{2-\epsilon}$ provided $k\geq 1$. Therefore,
$\sum_{n=a}^bf(n)$ and $\sum_{n=a}^bg(n)$ can be bounded above using integrals:
\begin{eqnarray*}
\frac{1}{2}\bigg|\sum_{n=1}^{\lfloor L^{2-\epsilon} \rfloor} E_{n}^{(1)}(L)
+\!\!\!\!\!\!\!\sum_{n=\lfloor L^{2-\epsilon} \rfloor+1}^\infty\!\!\!\!\!\!E_{n}^{(2)}(L)\bigg|
\leq \frac{C_\gamma}{2}\int_{L^{2-\epsilon}}^\infty\frac{dx}{x^{3/2-\gamma}}+
\sqrt{\frac{8L^2}{\pi}}\frac{1}{1-e^{-2L^\epsilon}} \int_0^{1+L^{2-\epsilon}}\!\!\!\!\frac{e^{-\frac{2L^2}{x}}}{x^{3/2}}dx 
\nonumber\\
\leq \frac{C_\gamma}{1-2\gamma}L^{-(2-\epsilon)\frac{1-2\gamma}{2}}
+\sqrt{\frac{8L^2}{\pi}}\frac{1}{1-e^{-2L^\epsilon}} 
\frac{1+L^{2-\epsilon}}{L^{3-3\epsilon/2}}e^{-\frac{2L^2}{1+L^{2-\epsilon}}}.
\end{eqnarray*}
Since $\gamma>0$ can be chosen to be arbitrarily small, we conclude from the above that
\bea
\frac{1}{2}\bigg|\sum_{n=1}^{\lfloor L^{2-\epsilon} \rfloor} E_{n}^{(1)}(L)
+\!\!\!\!\!\!\!\sum_{n=\lfloor L^{2-\epsilon} \rfloor+1}^\infty\!\!\!\!\!\!E_{n}^{(2)}(L)\bigg|=o\left(L^{-1+\epsilon/2+}\right).
\eea
Therefore,
\bea\label{rlterm3}
R(L)&=&\frac{1}{4\pi}\sum_{n=1}^{\lfloor L^{2-\epsilon}\rfloor} \frac{1}{n}\sum_{k=1}^{n-1} \frac{1}{\sqrt{k(n-k)}}
+\frac{1}{2}\sum_{n=\lfloor L^{2-\epsilon} \rfloor+1}^\infty\left(\frac{1}{2n}-\sqrt{\frac{2}{\pi n^3}}L \Omega 
\left(\frac{2L^2}{\pi n} \right) \right)\nonumber\\
&-&\frac{1}{4}\log\left(2L^2\right)
+o\left(L^{-1+\epsilon/2+}\right).
\eea 
Using integral bounds it is elementary to establish an estimate,
\bea
\sum_{k=1}^{n-1} \frac{1}{\sqrt{k(n-k)}}=\pi+O(n^{-1/2}),
\eea 
which can be used to re-write $R(L)$ as follows:
\bea\label{rlterm4}
R(L)&=&\frac{1}{4\pi}\sum_{n=1}^{\infty} \frac{1}{n}
\left(\sum_{k=1}^{n-1} \frac{1}{\sqrt{k(n-k)}}-\pi\right)
+\frac{1}{2}\sum_{n=\lfloor L^{2-\epsilon} \rfloor+1}^\infty\left(\frac{1}{2n}-\sqrt{\frac{2}{\pi n^3}}L \Omega 
\left(\frac{2L^2}{\pi n} \right) \right)\nonumber\\
&+&\frac{1}{4}\sum_{n=1}^{\lfloor L^{2-\epsilon}\rfloor} \frac{1}{n}
-\frac{1}{4}\log\left(2L^2\right)
+o\left(L^{-1+\epsilon/2+}\right).
\eea 
Recall a classical result for the sum of harmonic series \cite{wwt},
\bea\label{em}
\sum_{n=1}^N \frac{1}{n}=\log N+\gamma+O(N^{-1}),
\eea
where $\gamma$ is the Euler-Masceroni constant. Using (\ref{em}) in (\ref{rlterm4})
we find
\bea\label{rlterm5}
R(L)&=&\frac{1}{4\pi}\sum_{n=1}^{\infty} \frac{1}{n}
\left(\sum_{k=1}^{n-1} \frac{1}{\sqrt{k(n-k)}}-\pi\right)
+\frac{1}{2}\sum_{n=\lfloor L^{2-\epsilon} \rfloor+1}^\infty\left(\frac{1}{2n}-\sqrt{\frac{2}{\pi n^3}}L \Omega 
\left(\frac{2L^2}{\pi n} \right) \right)\nonumber\\
&-&\frac{\epsilon}{4}\log L+\gamma/4-\frac{1}{4}\log 2
+o\left(L^{-1+\epsilon/2+}\right).
\eea 
An application of the mean value theorem to the terms of the second sum on the right hand side
of (\ref{rlterm5}) leads to
\bea\label{penal}
&&\sum_{n=\lfloor L^{2-\epsilon} \rfloor+1}^\infty \!\!\left(\frac{1}{2n}-\sqrt{\frac{2}{\pi n^3}}L \Omega \left(\frac{2L^2}{\pi n} \right) \right)
=\int_{L^{-\epsilon} }^\infty \!\!\!\! dx
\left(\frac{1}{2x}-\sqrt{\frac{2}{\pi x^3}} \Omega
\left(\frac{2}{\pi x} \right) \right)+o(L^{-2+\epsilon})\nonumber\\
&&=\frac{\epsilon}{2}\log L
-\int_{0}^1  \!\!\!\!dx
\sqrt{\frac{2}{\pi x^3}} \Omega
\left(\frac{2}{\pi x} \right) +
\int_{1}^\infty \!\!\!\!\!\! dx
\left(\frac{1}{2x}-\sqrt{\frac{2}{\pi x^3}} \Omega
\left(\frac{2}{\pi x} \right) \right)
+o(L^{-2+\epsilon}).
\eea
Substituting (\ref{penal}) into the r.h.s of (\ref{rlterm5}) we discover that
\bea\label{rlterm6}
R(L)\!&=&\!\gamma/4\!-\!\frac{1}{4}\log 2
+\frac{1}{4\pi}\sum_{n=1}^{\infty} \frac{1}{n}
\left(\sum_{k=1}^{n-1} \frac{1}{\sqrt{k(n-k)}}-\pi\right)\nonumber\\
\!\!&-&\!\!\frac{1}{2}\int_{0}^1  \!\!\!\!dx
\sqrt{\frac{2}{\pi x^3}} \Omega
\left(\frac{2}{\pi x} \right) \!\!+\!\!
\frac{1}{2}\int_{1}^\infty \!\!\!\!\!\! dx
\left(\frac{1}{2x}\!-\!\sqrt{\frac{2}{\pi x^3}} \Omega
\left(\frac{2}{\pi x} \right) \right)
\!\!+\!\!o\left(\!L^{-1+\epsilon/2+}\!\right)\!\!.
\eea 
As $\epsilon>0$ is arbitrary, we conclude that the magnitude of the error term is $o(L^{-1+})$.

In principle, (\ref{rlterm6}) gives an answer for the $O(1)$ term in the expansion
of $\pr(\lm<-L)$. It can however be considerably simplified, which probably means
that the calculation detailed above can also be significantly streamlined. The rest of the
proof is an exact calculation based on the relation between the Euler-Masceroni constant
and Jacobi's theta functions.   

The calculation is based on the following two remarks: firstly,
\bea\label{emtheta}
\gamma=\log(4\pi)-2+2 \int_1^\infty (1+\sqrt{t}) \frac{\Omega(t)}{t}dt,
\eea
see \cite{choi} containing this as well as a large collection of other expressions
for the Euler-Masceroni constant. Formula (\ref{emtheta}) follows from combining a
more
standard expression for $\gamma$ in terms of $\zeta$-function,
\[
\gamma=\lim_{s\rightarrow 1} \left[\zeta(s)-\frac{1}{s-1}\right],
\]  
see \cite{wwt} for the derivation, and Riemann's integral representation of $\zeta$,
\[
\zeta (s)=\frac{\pi^{s/2}}{s(s-1)\Gamma(s/2)}+\frac{\pi^{s/2}}{\Gamma(s/2)}\int_1^\infty
\left(t^{(1-s)/2}+t^{s/2}\right) \frac{\Omega(t)}{t} dt,
\]
see e.g. \cite{edwards}. Secondly, 
\bea\label{mod}
1+2\Omega\left( t^{-1}\right)=\sqrt{t} \left(1+2\Omega(t)\right),
\eea
which follows from the standard transformation properties of the theta function
and can be proved directly using Poisson summation formula, see \cite{wwt} for review.
Now let us modify the r.h.s. of (\ref{rlterm6}) as follows: express $\gamma$ using (\ref{emtheta}), 
change variables in the penultimate integral according to $t=\frac{2}{\pi x}$, apply (\ref{mod})
 to the integrand of the last integral. The result is
 \bea\label{rlterm7}
R(L)+o(L^{-1+})\!\!\!\!\!&=&\!\frac{1}{4}\log 2\pi-\frac{1}{2}+\sqrt{\frac{1}{2\pi}}
+\frac{1}{4\pi}\sum_{n=1}^{\infty} \frac{1}{n}
\left(\sum_{k=1}^{n-1} \frac{1}{\sqrt{k(n-k)}}-\pi\right)\nonumber\\
\!\!&-&\!\!\frac{1}{2}\int_{2/\pi}^1  \!\!\!\!dt
\frac{\Omega(t)}{\sqrt{t}}
\!\!+\!\!
\frac{1}{2}\int_{1}^{\pi/2} \!\!\!\!\!\! dt
\frac{\Omega(t)}{t}
=\!\frac{1}{4}\log 2\pi-\frac{1}{2}+\sqrt{\frac{1}{2\pi}}\nonumber\\
&+&\frac{1}{4\pi}\sum_{n=1}^{\infty} \frac{1}{n}
\left(\sum_{k=1}^{n-1} \frac{1}{\sqrt{k(n-k)}}-\pi\right)
+\frac{1}{2}\int_{1}^{\pi/2} \!\!\!\!\!\! dt
\frac{\sqrt{t}\Omega(t)-\Omega(t^{-1})}{t^{3/2}}\nonumber\\
&=&\!\frac{1}{2}\log 2
+\frac{1}{4\pi}\sum_{n=1}^{\infty} \frac{1}{n}
\left(\sum_{k=1}^{n-1} \frac{1}{\sqrt{k(n-k)}}-\pi\right),
\eea 
where the last equality follows from another application of (\ref{mod})
to the integral in the previous expression.
Theorem \ref{thm} is proved.\qed
\section{The proof of probabilistic lemmas}\label{proof_lemmas}
\subsection{Lemma \ref{thm:cyclic}}
Let $(X_k)_{1\leq k \leq 2n}$ for $n=1,2,\ldots$, be a sequence of independent indentically
distributed random variables with a continuous density.
Let $S=(S_k)_{0\leq k \geq 2n}$ be the associated random walk started at zero,
\bea
S_{k}=\sum_{m=1}^k X_m.
\eea
We will frequently consider the walk conditioned to be at $0$ at time $2n$, so that $S_{2n}=0$.
Let $S^{(p)}$ be the random walk associated with the cyclic shift of the increments $X$'s
by $p$ steps to the right,
\bea
S^{(p)}_k=\sum_{m=1}^k X_{m+p},
\eea
where the addition in the time indices is performed modulo $2n$. As it is easy to check,
under the conditioning that $S_{2n}=0$,
\bea\label{csh}
S_{k}^{(p)}=S_{k+p}-S_{p},~S_{2n}^{(p)}=0.
\eea
Under the conditioning $S_{2n}=0$, for any  $p$, $(S_{k}^{(p)})_{0\leq k \leq 2n}$ is a bridge, whose law is $p$-independent,
\bea\label{ci}
(S_{k}^{(p)})_{0\leq k \leq 2n}\stackrel{(d)}{\sim} (S_{k}^{(q)})_{0\leq k \leq 2n},~0\leq p,q \leq 2n-1.
\eea
Recall that $\tau_0$ be the first even time the value of the bridge $S$ becomes negative.
We write $\tau_0^{(p)}$ for the corresponding exit time for the bridge $S^{(p)}$. Similarly,
$M_{2n}^{(p)}$ is the odd time maximum of the walk $S^{(p)}$.
Using the fact that the even time global minimum of the bridge is unique almost surely, one can easily show
that
\bea\label{indsum}
\sum_{p=0}^{n-1} \ind\left(\tau_{0}^{(2p)}=2n\right)=1~a.~s.
\eea
Notice that the event $\tau_0^{(2p)}=2n$ corresponds to the shift to the time at which the even time global
minimum of the random walk has been achieved.
Let 
\bea\label{minim}
m_{2n}=\min_{0\leq k\leq n} S_{2k}.
\eea
Therefore,
\bea
p_{n}(L)&:=&\E\left(\min(L,M_{2n}) \ind(\tau_0=2n)\delta_0( S_{2n}) \right)
\nonumber
\\
&\stackrel{(\ref{ci})}{=}&
\frac{1}{n}\sum_{p=0}^{n-1} 
\E\left(\min(L,M_{2n}^{(2p)}) \ind(\tau^{(2p)}_0=2n)\mid S_{2n}^{(2p)}=0 \right)\Pr(S_{2n}^{(2p)}\in d0)
\nonumber\\
&\stackrel{(\ref{csh})}{=}&\frac{1}{n}\sum_{p=0}^{n-1} 
\E\left(\min(L,M_{2n}-S_{2p}) \ind(\tau^{(2p)}_0=2n)\mid S_{2n}=0 \right)\Pr(S_{2n}\in d0)
\nonumber\\
&\stackrel{(\ref{minim})}{=}&\frac{1}{n}\sum_{p=0}^{n-1} 
\E\left(\min(L,M_{2n}-m_{2n}) \ind(\tau^{(2p)}_0=2n)\mid S_{2n}=0 \right)\Pr(S_{2n}\in d0)
\nonumber\\
&=& \frac{1}{n}
\E\left(\min(L,M_{2n}-m_{2n}) 
\sum_{p=0}^{n-1}\ind(\tau^{(2p)}_0=2n)\mid S_{2n}=0 \right)\Pr(S_{2n}\in d0)
\nonumber\\
&\stackrel{(\ref{indsum})}{=}& 
\frac{1}{n}\E\left(\min(L,M_{2n}-m_{2n}) 
\mid S_{2n}=0 \right)\Pr(S_{2n}\in d0)
\nonumber\\
&=&\frac{1}{n}\E\left(\min(L,M_{2n}-m_{2n}) 
\delta_0(S_{2n})\right).
\nonumber
\eea
 \qed\\
 \\
 \begin{remark}\label{leadtrm} 
 Using the notations developed for the proof, it is very easy to rederive (\ref{eqlt}),
 even though the  computation presented below uses no new ideas in comparison with \cite{pop}:
 \begin{eqnarray*}
 \E\left(\delta_0(S_{\tau_0})\right)&=&\sum_{n=1}^\infty \E\left(\ind(\tau_0=2n)\delta_0(S_{2n})\right)
 =
 \sum_{n=1}^\infty \E\left(\ind(\tau_0=2n)\mid S_{2n}=0)\Pr(S_{2n}\in d0\right)
\\
&=& \sum_{n=1}^\infty \frac{1}{n}\sum_{p=0}^{n-1}
\E\left(\ind(\tau^{(2p)}_0=2n)\mid S^{(2p)}_{2n}=0\right)\Pr(S_{2n}\in d0)
\\
&=&\sum_{n=1}^\infty \frac{1}{n}\sum_{p=0}^{n-1}\E\left(\ind(\tau^{(2p)}_0=2n
\mid S_{2n}=0\right)\Pr(S_{2n}\in d0)\\
&=&
\sum_{n=1}^\infty \frac{1}{n}\E\left(\sum_{p=0}^{n-1}\ind(\tau^{(2p)}_0=2n)
\mid S_{2n}=0\right)\Pr(S_{2n}\in d0)
\\
&\stackrel{(\ref{indsum})}{=}&\sum_{n=1}^\infty \frac{1}{n}\E\left(1
\mid S_{2n}=0\right)\Pr(S_{2n}\in d0)
=\sum_{n=1}^\infty \frac{1}{n}\Pr(S_{2n}\in d0)\\
&=&\sum_{n=1}^\infty \frac{1}{n}\frac{1}{\sqrt{2\pi n}}
=\frac{\zeta(3/2)}{\sqrt{2\pi}}.
 \end{eqnarray*} 
 \end{remark}
 \subsection{Lemma \ref{errbounds}}
 \subsubsection{$n\leq L^{2-\epsilon}$} 
 Starting from formula (\ref{pncyc}) of Lemma \ref{thm:cyclic}, 
 \bea
 p_n(L)=\frac{1}{n}\E\left(\left(M_{2n}-m_{2n}\right)\delta_0(B_{2n})\right)+E^{(1)}_n(L),
 \eea
 where
 \bea
 E^{(1)}_n(L)=-\frac{1}{n}\E\left(\left(M_{2n}-m_{2n}-L\right)_{+}\delta_{0}(B_{2n})\right),
 \eea
 and $x_{+}:=x\ind(x\geq 0)$. 
 
 We start with estimating the error term $E^{(1)}_n(L)$. Let $(B_t)_{t\geq 0}$ be the
 rate-$1/2$ Brownian motion. It follows from the definition of $M_{2n}$, $m_{2n}$ that
 \bea
 M_{2n} \leq \sup_{0\leq t\leq 2n} B_t,~ m_{2n}\geq \inf_{0\leq t \leq 2n} B_t.
 \eea 
 As the function $x\mapsto x_{+}$ is increasing,
 \bea\label{est1}
 |E^{(1)}_n(L)|&=&\frac{1}{n}\E\left(\left(M_{2n}-m_{2n}-L\right)_{+}\delta_{0}(B_{2n})\right)
 \nonumber\\
 &\leq& 
 \frac{1}{n}\E\left(\left(\sup_{0\leq t\leq 2n} B_t- \inf_{0\leq t \leq 2n} B_t-L\right)_{+}
 \delta_{0}(B_{2n})\right)
 \nonumber\\
 &\leq& 
 \frac{1}{n}\E\left(\left(\sup_{0\leq t\leq n} W_t- \inf_{0\leq t \leq n} W_t-L\right)_{+}\delta_{0}(W_{n})\right),
 \eea
 where $(W_t)_{t\geq 0}$ is the standard Brownian motion. Thus the error term $\E^{(1)}_n(L)$
 is bounded by an expectation w. r. t. to the Wiener measure, which can be computed using
 the known joint distribution of the supremum, infimum and the final position of the Brownian motion,
 \bea\label{trivar}
 \pr\left(
 \inf_{0\leq t \leq n} W_t \geq a
 \sup_{0\leq t\leq n}\leq b, W_t \in d0 \right)=
 \sum_{k\in \Z}\frac{1}{\sqrt{2\pi t}}\left(e^{-\frac{2k^2}{t}(b-a)^2}-e^{\frac{2}{t}(b-k(b-a))^2} \right),
 \eea
 where $a\leq 0, b \geq 0$,
 see e.g. \cite{billingsey}. Applying (\ref{trivar}) to the computation of the r. h. s. of (\ref{est1}) 
 one finds after some work that
\bea
 |E^{(1)}_n(L)|\leq \frac{1}{\sqrt{8\pi nL^2}}e^{-\frac{2L^2}{n}}+\sqrt{\frac{2}{\pi n^3}}L
 \sum_{k=1}^\infty e^{-\frac{2k^2L^2}{n}}.
\eea
Finally, approximating $e^{-\frac{2k^2L^2}{n}}\leq e^{-\frac{2kL^2}{n}}$ for $k\geq 1$, summing
the resulting geometric series and combining the terms leads to
formula (\ref{e1}) of Lemma \ref{errbounds}.

To finish the proof of (\ref{plow}), we need to calculate $\frac{1}{n}\E\left(\left(M_{2n}-m_{2n}\right)\delta_0(B_{2n})\right)$.
Using reflection symmetry,
\bea\label{mc11}
\E\left(m_{2n}\delta_0(B_{2n})\right)=
-\E\left(\max_{0\leq k\leq n} (B_{2k})\delta_0(B_{2n})\right)
=-\E\left(\max_{0\leq k\leq n} (W_{k})\delta_0(W_{n})\right),
\eea
where $(W_k)_{k \geq 0}$ is the random walk with $N(0,1)$ increments.
Let us denote the increments of random walk $B$ by $X_1, X_2, \ldots \sim N(0,1/2)$,
let $Y_1, Y_2,, \ldots \sim N(0,1)$ be the increments of the walk $W$. Then
\bea\label{mc10}
\!\!\!\!\!\!\!\!\!&&\E\left(\left(M_{2n}\right)\delta_0(B_{2n})\right)=
\E\left(\max_{0\leq k\leq n-1} (B_{2k+1})\delta_0(B_{2n})\right)\nonumber\\
&=&\!\!\!\!\E\left(\max(X_1, X_1+X_2+X_3, \ldots, X_1+X_2+\ldots +X_{2n-1} )\delta_0(B_{2n})\right)\nonumber\\
&=&\!\!\!\!\E\left(X_1+\max(0, X_2+X_3, X_2+X_3+X_4+X_5, \ldots, X_3+X_5+\ldots +X_{2n-1} )\delta_0(B_{2n})\right)\nonumber\\
&=&\!\!\!\!\E\left(\max(0, Y_1, Y_1\!+\!Y_2, \ldots, Y_1\!+\!Y_2\!+\!\ldots \!+\!Y_{n-1} )\delta_0(W_{n})\right)
\!\!=\!\!\E\left(\max_{0\leq k \leq n} (W_k) \delta_0(W_n)\right).
\eea
Substituting (\ref{mc11}) and (\ref{mc10}) into $\frac{1}{n}\E\left(\left(M_{2n}-m_{2n}\right)\delta_0(B_{2n})\right)$
we find
\bea
\frac{1}{n}\E\left(\left(M_{2n}-m_{2n}\right)\delta_0(B_{2n})\right)=\frac{2}{n}\E(\max_{0\leq k\leq n} W_{k}\delta_0
(W_n))=
\frac{1}{2\pi n}\sum_{k=1}^{n-1} \frac{1}{\sqrt{k(n-k)}},
\eea
where the last step used formula (1.11a) from \cite{kac} (Lemma \ref{lemma_kac} of the present paper). Formula (\ref{plow}) of Lemma \ref{errbounds} is proved.\\
 \subsubsection{$n \geq L^{2-\epsilon}$}
 For large values of $n$ it is natural to approximate the Gaussian random walk with Brownian motion
 to re-write (\ref{pncyc}) as follows:
 \bea\label{spoint}
 p_n(L)=\frac{1}{n} \E\left( \min\left(L,\sup_{0\leq t \leq 2n}B_t-\inf_{0\leq t\leq 2n}B_t\right)
 \delta_0\left(B_{2n}\right)\right)+E_n^{(2)}(L),
 \eea
 where the correction term is
 \bea\label{est5}
 &&E_n^{(2)}(L)\!\\
 &=&\!\frac{1}{n}
 \E\left(\left(\min\left(L,\!\!\!\!\sup_{0\leq 2k+1 \leq 2n}\!\!\!\!B_{2k+1}-\!\!\inf_{0\leq 2k\leq 2n}
 \!\!\!\!B_{2k}\right)\!\!-\!\!
 \min\left(L,\!\!\sup_{0\leq t \leq 2n}\!\!\!\!B_t-\!\!\inf_{0\leq t\leq 2n}\!\!\!\!B_t\right)\right)\!\!\delta_0\left(B_{2n}\right)\right).
 \nonumber
 \eea
 For any $L,x,y \in \R$, 
 \[
 |\min(L,x)-\min(L,y)|\leq |x-y|,
 \]
 which allows us to bound (\ref{est5}) as below:
 \bea\label{est26}
 E_n^{(2)}(L) &\leq& \frac{1}{n}\E\left(\bigg|\sup_{0\leq t \leq 2n}B_t-\sup_{0\leq k< n}
 B_{2k+1}\bigg| \delta_0\left(B_{2n}\right)\right)\nonumber \\
 &+&
 \frac{1}{n}\E\left(\bigg|\inf_{0\leq t \leq 2n}B_t-\inf_{0\leq k\leq n}
 B_{2k}\bigg| \delta_0\left(B_{2n}\right)\right).
 \eea
 The two terms on the r.h.s. are very similar and can be bounded by the same function of
 the index $n$. We will present the derivation of the bound for the first term only.
 In what follows, $(W_t)_{t \geq 0}$ is the standard Brownian motion, $(WB_t)_{0\leq t \leq 1}$
 is the Brownian bridge. 
 
 Rescaling time,
 \bea\label{est13}
&&\frac{1}{n}\E\left(\bigg|\sup_{0\leq t \leq 2n}B_t-\!\!\!\!\sup_{0\leq 2k+1\leq 2n}\!\!\!\!\!
 B_{2k+1}\bigg| \delta_0\left(B_{2n}\right)\right)\!\!
 =\!\!
 \frac{1}{n}\E\left(\bigg|\sup_{t \in [0,1]}\!\!W_t-\!\!\!\!\!\!\sup_{t \in \{\frac{k}{n}
 +\frac{1}{2n}\}_{k=0}^{n-1}}\!\!\!\!\!\!\!\!
 W_t\bigg| \delta_0\left(W_{1}\right)\right)
 \nonumber\\
 &=&\!\!\!\!
 \sqrt{\frac{1}{2\pi n^2}}
 \E\left(\bigg|\sup_{t \in [0,1]}\!\!W_t-\!\!\!\!\sup_{t \in \{\frac{k}{n}+
 \frac{1}{2n}\}_{k=0}^{n-1}}\!\!\!\!\!\!\!\!
 W_t\bigg| \bigg| W_{1}=0\right)
 \nonumber\\
  &=&\!\!
 \sqrt{\frac{1}{2\pi n^2}}
 \E\left(\bigg|\sup_{t \in [0,1]}\!\!WB_t-\!\!\!\!\sup_{t \in \{\frac{k}{n}+
 \frac{1}{2n}\}_{k=0}^{n-1}}\!\!\!\!\!\!\!\!
 WB_t\bigg|\right).
 \eea
 Recall the following fact about Brownian motions hence the Brownian bridges: 
 for any fixed $\gamma>0$ there is a non-negative
 random variable $H_\gamma$ defined 
 on the same probability space as the bridge itself such that, 
 \bea\label{estbb}
 |WB_{t}-WB_{\tau}|\leq H_\gamma |t-\tau|^{\frac{1}{2}-\gamma}, \mbox{ for all } t,\tau \in [0,1].
 \eea
 Moreover, $\E(H_\gamma)<\infty$, see e.g. \cite{billingsey}. Exploiting (\ref{estbb}) to bound  the r.h.s. of (\ref{est13})
 one finds that
 \bea\label{est27}
 \frac{1}{n}\E\left(\bigg|\sup_{0\leq t \leq 2n}B_t-\!\!\!\!\sup_{0\leq 2k+1\leq 2n}\!\!\!\!\!
 B_{2k+1}\bigg| \delta_0\left(B_{2n}\right)\right)\!\! \leq 
 \frac{\E(H_\gamma)}{\sqrt{2\pi n^2}} \left(\frac{1}{2n} \right)^{\frac{1}{2}-\gamma}.
 \eea
 The second second term on the r.h.s. of (\ref{est26}) obeys the same bound. Thus combining
 (\ref{est26}) with (\ref{est27}) we conclude that
 \bea
 E_{n}^{(2)}(L)\leq C_\gamma n^{-3/2+\gamma},
 \eea
 where $C_\gamma=2\frac{\E(H_\gamma)}{\sqrt{2\pi }} \left(\frac{1}{2} \right)^{\frac{1}{2}-\gamma}$
 is an $n$-independent constant. The estimate (\ref{e2})
 of Lemma \ref{errbounds} is proved.
 
In order to calculate the leading term in the expression (\ref{spoint}) for $p_{n}(L)$ we just need to 
evaluate the expectation of $ \min\left(L,\sup_{0\leq t \leq 2n}B_t-\inf_{0\leq t\leq 2n}B_t\right)$
using the distribution (\ref{trivar}). The answer is
\bea
\frac{1}{n} \E\left( \min\left(L,\sup_{0\leq t \leq 2n}B_t-\inf_{0\leq t\leq 2n}B_t\right)
 \delta_0\left(B_{2n}\right)\right)=\frac{1}{2n}-\frac{2L}{\sqrt{2\pi n^3}}\sum_{k=1}^\infty e^{-\frac{2k^2}{n}L^2}.
\eea
Formula (\ref{phi}) of Lemma \ref{errbounds} is proved.
\qed
 \subsection{Lemma \ref{lemma_kac}}
 Fix $k<n$ and define 
 \bea
 E_k:=\E\left(\max(0, S_1, \ldots, S_k) \mid S_n=0\right).
 \eea
 Decomposing the expectation according to the events $S_k< 0$ and $S_k>0$ we find
 \begin{eqnarray*}
 E_k&=&\E\left(\ind(S_k>0)\max(S_1, \ldots, S_k) \mid S_n=0\right)+
 \E\left(\ind(S_k<0)\max(0, S_1, \ldots, S_k) \mid S_n=0\right)\\
 &=&\E\left(\ind(S_k>0)X_1 \mid S_n=0\right)+(
 \E(\ind(S_k>0)\max(0,X_2,X_2+X_3 \ldots, S_k-X_1) \mid S_n=0)\\
 &+&
 \E(\ind(S_k<0)\max(0, S_1, \ldots, S_{k-1}) \mid S_n=0) )
 =\frac{1}{k}\E(\ind(S_k>0)S_k\mid S_n=0)+E_{k-1},
 \end{eqnarray*}
 where the last equality is due to the invariance of the law of the increments $X_1, X_2, \ldots X_n$
 with respect to a permutation of the first $k$ increments
 and the $X\rightarrow -X$ symmetry of the distribution of increments. 
 Solving the resulting difference equations
 for $E_k$'s with the initial condiiton $E_0=0$, we find
 \bea\label{dyson}
 \E\left(\max(0, S_1, \ldots, S_{n-1}) \mid S_n=0\right)=\sum_{k=1}^{n-1}\frac{1}{k}\E(\ind(S_k>0)S_k\mid S_n=0).
 \eea
 (Formula (\ref{dyson}) is attributed in \cite{kac} to Freeman Dyson.) Therefore,
 \bea
 \rho^{(n)}(0)\E(\max(0, S_1, &\ldots&, S_{n-1}) \mid S_n=0)\stackrel{(\ref{dyson})}{=}
 \sum_{k=1}^{n-1}\frac{1}{k}\E(\ind(S_k>0)S_k\mid S_n=0)\rho^{(n)}(0)\nonumber\\
 &=&\sum_{k=1}^{n-1}\frac{1}{k}\E((S_k)_{+}\delta_0( S_n))
=\sum_{k=1}^{n-1}\frac{1}{k}\int_{0}^\infty dx x \rho^{(k)}(x)\rho^{(n-k)}(x)\nonumber\\
&=& \frac{n}{2}\sum_{k=1}^{n-1}\frac{1}{k(n-k)}\int_{0}^\infty dx x \rho^{(k)}(x)\rho^{(n-k)}(x).
\eea
 Lemma \ref{lemma_kac} is proved. \qed \\
\\
{\bf Acknowledgement.} We are grateful to Thomas Bothner for illuminating discussions.

\bibliographystyle{amsplain}

\end{document}